\documentclass[12pt]{article}
\input epsf

\title{On the number of extremal surfaces.}
\author{\sc Alina Vdovina}
\date{
\small {Mathematisches Institut} \\
\small {Beringstr. 1, Bonn 53115} \\
\small {e-mail: alina@math.uni-bonn.de}
}

\begin{document}
\maketitle

\begin{abstract}

Let $X$ be a compact Riemann surface of genus $\geq 2$
of constant negative curvature $-1$. An extremal disk is
an embedded (resp. covering) disk of maximal (resp.
minimal) radius. A surface containing an extremal disk
is an {\em extremal surface}.

This paper gives formulas enumerating extremal surfaces
of genus $\geq 4$ up to isometry. We show also that the
isometry group of an extremal surface is always cyclic
of order 1, 2, 3 or 6.

\end{abstract}

\section*{Introduction}

Let $X$ be a compact Riemann surface of genus $\geq 2$
of constant negative curvature $-1$.
We consider the maximal radius of an embedded metric
disk in $X$ and the minimal radius of a disk covering $X$.

An extremal disk is an embedded (resp. covering) disk of maximal
(resp. minimal) radius. A surface containing an extremal
disk is an {\em extremal surface}.
C.Bavard \cite{[Ba]} proved, that if a surface contains
an embedded disk of maximal radius if and only if it contains
a covering disk of minimal radius and that extremal
surfaces are modular surfaces.

The radius $R_g$ of an extremal
embeddes disk, as well as the radius 
$C_g$ of an extremal covering disk were computed in  \cite{[Ba]}:

$$R_g=cosh^{-1}(1/2sin \beta_g), \beta_g=\pi/{(12g-6)}.$$
 $$C_g=cosh^{-1}(1/\sqrt3 tan \beta_g), \beta_g=\pi/{(12g-6)}.$$
C.Bavard \cite{[Ba]} also proved, that discs of maximal radius occur
in those surfaces which admit as Dirichlet domain 
a regular polygon with the largest possible number of
sides $12g-6$ and extremal surfaces occur in every genus.

We give an explicit construction of all extremal surfaces
of genus $g\geq 4$. 
We show, that for genus $g \geq 4$ the isometry group
of an extremal surface is always cyclic of order 1,2,3 or 6
and we give an explicit formula for the number
of nonisometric extremal surfaces.
Also we give explicit formulas of extremal surfaces
having exactly $d$ automorphisms, where $d$ is 1,2,3 or 6.
Form those formulas one can see that asymptotically
almost all extremal surfaces have no isometries.
In particular, for $d=1$ we have a big family of
explicitly constructed surfaces with no automorphisms.
Let's note, that another families of surfaces with no
automorphisms were considered in \cite{[E]}, \cite{[T]}.
The questions of explicit construction, enumeration and
describtion of isometries of genus $2$ extremal
surfaces were solved in \cite{[GG3]}.
 We will show, that oriented maximal
 Wicks forms
and extremal surfaces
are in bijection for $g \geq 4$. For $g=2$ this bijection was
proved in \cite{[GG1]} and for $g=3$ the question is still open.

 Section \ref{sec:1} formulates our main results and 
introduces oriented Wicks forms (cellular
 decompositions with only one face of oriented surfaces), our main tool.
 Wicks forms are canonical forms for products of commutators
 in free groups \cite{[V]}. 

 Section \ref{sec:2} contains a few facts concerning oriented maximal
 Wicks forms.

 Section \ref{sec:3} contains the proof of our main results.

\section{Main results}
\label{sec:1}


 {\bf Definition 1.1.}
 An {\it oriented Wicks form\/} is a cyclic word $w= w_1w_2\dots w_{2l}$
 (a cyclic word is the orbit of a linear word under cyclic permutations)
 in some alphabet $a_1^{\pm 1},a_2^{\pm 1},\dots$ of letters
 $a_1,a_2,\dots$ and their inverses $a_1^{-1},a_2^{-1},\dots$ such that
\begin{itemize}
\item[(i)] if $a_i^\epsilon$ appears in $w$ (for $\epsilon\in\{\pm 1\}$)
 then $a_i^{-\epsilon}$ appears exactly once in $w$,
\item[(ii)] the word $w$ contains no cyclic factor (subword of
 cyclically consecutive letters in $w$) of the form $a_i a_i^{-1}$ or
 $a_i^{-1}a_i$ (no cancellation),
\item[(iii)] if $a_i^\epsilon a_j^\delta$ is a cyclic factor of $w$ then
 $a_j^{-\delta}a_i^{-\epsilon}$ is not a cyclic factor of $w$ 
(substitutions of the form
 $a_i^\epsilon a_j^\delta\longmapsto x,
 \quad a_j^{-\delta}a_i^{-\epsilon}\longmapsto x^{-1}$ are impossible).
\end{itemize}

 An oriented Wicks form $w=w_1w_2\dots$ in an alphabet $A$ is
 {\em isomorphic\/} to $w'=w'_1w'_2$ in an alphabet $A'$ if
 there exists a bijection $\varphi:A\longrightarrow A'$ with
 $\varphi(a^{-1})=\varphi(a)^{-1}$ such that $w'$ and
 $\varphi(w)=\varphi(w_1)\varphi(w_2)\dots$ define the
 same cyclic word.

 An oriented Wicks form $w$ is an element of the commutator subgroup
 when considered as an element in the free group $G$ generated by
 $a_1,a_2,\dots$. We define the {\em algebraic genus\/} $g_a(w)$ of
 $w$ as the least positive integer $g_a$ such that $w$ is a product
 of $g_a$ commutators in $G$.

 The {\em topological genus\/} $g_t(w)$ of an oriented Wicks
 form $w=w_1\dots w_{2e-1}w_{2e}$ is defined as the topological
 genus of the oriented compact connected surface obtained by
 labeling and orienting the edges of a $2e-$gone (which we
 consider as a subset of the oriented plane) according to
 $w$ and by identifying the edges in the obvious way.

{\bf Proposition 1.1.}
{\sl The algebraic and the topological genus of an oriented Wicks
 form coincide (cf. \cite{[C],[CE]}).}

 We define the {\em genus\/} $g(w)$ of an oriented
 Wicks form $w$ by $g(w)=g_a(w)=g_t(w)$.

 Consider the oriented compact surface $S$ associated to an oriented 
 Wicks form $w=w_1\dots w_{2e}$. This surface carries an immerged graph
 $\Gamma\subset S$ such that $S\setminus \Gamma$ is an open polygon
 with $2e$ sides (and hence connected and simply connected).
 Moreover, conditions (ii) and (iii) on Wicks form imply that $\Gamma$ 
 contains no vertices of degree $1$ or $2$ (or equivalently that the
 dual graph of $\Gamma\subset S$ contains no faces which are $1-$gones
 or $2-$gones). This construction works also
 in the opposite direction: Given a graph $\Gamma\subset S$
 with $e$ edges on an oriented compact connected surface $S$ of genus $g$
 such that $S\setminus \Gamma$ is connected and simply connected, we get
 an oriented Wicks form of genus $g$ and length $2e$ by labeling and 
 orienting the edges of $\Gamma$ and by cutting $S$ open along the graph
 $\Gamma$. The associated oriented Wicks form is defined as the word
 which appears in this way on the boundary of the resulting polygon
 with $2e$ sides. We identify henceforth oriented Wicks
 forms with the associated immerged graphs $\Gamma\subset S$,
 speaking of vertices and edges of oriented Wicks form.

 The formula for the Euler characteristic
 $$\chi(S)=2-2g=v-e+1$$
 (where $v$ denotes the number of vertices and $e$ the number
 of edges in $\Gamma\subset S$) shows that
 an oriented Wicks form of genus $g$ has at least length $4g$
 (the associated graph has then a unique vertex of degree $4g$
 and $2g$ edges) and at most length $6(2g-1)$ (the associated
 graph has then $2(2g-1)$ vertices of degree three and
 $3(2g-1)$ edges).

 We call an oriented Wicks form of genus $g$ {\em maximal\/} if it has
 length $6(2g-1)$. Oriented maximal Wicks forms are dual to 1-vertex
 triangulations. This can be seen by cutting the oriented surface $S$ 
 along $\Gamma$, hence obtaining a polygon $P$ with $2e$ sides. 
 We draw a star $T$ on $P$ which joins an interior point of $P$
 with the midpoints of all its sides. Regluing $P$ we recover $S$
 which carries now a 1-vertex triangulation given by $T$ and each
 1-vertex triangulation is of this form for some oriented maximal
 Wicks form (the immerged graphs $T\subset S$ and $\Gamma\subset S$
 are dual to each other: faces of $T$ correspond to vertices of
 $\Gamma$ and vice-versa. Two faces of $T$ share a common edge if
 and only if the corresponding vertices of $\Gamma$ are adjacent).
 This construction shows that we can work indifferently with
 1-vertex triangulations or with oriented maximal Wicks forms.

 Similarly, cellular decompositions of oriented surfaces with
 one vertex and one face correspond to oriented minimal 
 Wicks forms and were enumerated in \cite{[CM]}. The dual of an
 oriented minimal Wicks form is again a (generally non-equivalent)
 oriented minimal Wicks form and taking duals yields hence an
 involution on the set of oriented minimal Wicks forms.

 A vertex $V$ (with oriented edges $a,b,c$ pointing toward $V$) is
 {\em positive\/} if
$$w=ab^{-1}\dots bc^{-1}\dots ca^{-1}\dots \quad {\rm or }\quad
 w=ac^{-1}\dots cb^{-1}\dots ba^{-1}\dots $$
 and $V$ is {\em negative\/} if    
 $$w=ab^{-1}\dots ca^{-1}\dots bc^{-1}\dots \quad {\rm or }\quad 
 w=ac^{-1}\dots ba^{-1}\dots ab^{-1}\dots \quad
 ..$$

 The {\em automorphism group\/} ${\rm Aut}(w)$ of an oriented
 Wicks form $$w=w_1w_2\dots w_{2e}$$ of length $2e$ is the group of all
 cyclic permutations $\mu$ of the linear word $w_1w_2\dots w_{2e}$ such
 that $w$ and $\mu(w)$ are isomorphic linear words (i.e. $\mu(w)$ is
 obtained from $w$ by permuting the letters of the alphabet). The group
 ${\rm Aut}(w)$ is a subgroup of the cyclic group ${\bf Z}/2e{\bf Z}$
 acting by cyclic permutations on linear words representing $w$.

 The automorphism group ${\rm Aut}(w)$ of an oriented Wicks
 form can of course also be described in terms of permutations on the
 oriented edge set induced by orientation-preserving
 homeomorphisms of $S$ leaving $\Gamma$ invariant. In particular
 an oriented maximal Wicks form and the associated dual
 1-vertex triangulation have isomorphic automorphism groups.

 We define the {\em mass\/} $m(W)$ of a finite set $W$ of oriented
 Wicks forms by
 $$m(W)=\sum_{w\in W}{1\over \vert{\rm Aut}(w)\vert}\quad .$$

 Let us introduce the sets
\par $W_1^g$: all oriented maximal
 Wicks forms of genus $g$ (up to equivalence),
\par $W^g_2(r)\subset W_1^g$: all oriented 
 maximal Wicks forms having an automorphism of order $2$ leaving 
 exactly $r$ edges of $w$ invariant by reversing their
 orientation. (This automorphism is the half-turn with respect to the
 \lq\lq midpoints" of these edges and exchanges the two adjacent
 vertices of an invariant edge.)
\par $W^g_3(s,t)\subset W_1^g$: all oriented maximal 
 Wicks forms having an automorphism of order $3$ leaving
 exactly $s$ positive and $t$ negative vertices invariant
 (this automorphism permutes cyclically the edges around
 an invariant vertex).
\par $W^g_6(3r;2s,2t)=W^g_2(3r)\cap W^g_3(2s,2t)$: all oriented
 maximal Wicks forms having an automorphism $\gamma$ of order 6
 with $\gamma^3$ leaving $3r$ edges invariant and $\gamma^2$
 leaving $2s$ positive and $2t$ negative vertices invariant
 (it is useless to consider the set $W_6^g(r';s',t')$ defined
 analogously since $3$ divides $r'$ and $2$ divides $s',t'$ if
 $W_6^g(r';s',t')\not=\emptyset$).

\par We define now the {\em masses\/} of these sets as
  $$\matrix
   {m_1^g\hfill&=&\displaystyle \sum_{w\in W_1^g}
    {1\over \vert{\rm Aut}(w)\vert}\quad ,\hfill \cr
    m_2^g(r)\hfill&=&\displaystyle \sum_{w\in W_2^g(r)}
    {1\over \vert{\rm Aut}(w)\vert}\quad ,\hfill \cr
    m_3^g(s,t)\hfill&=&\displaystyle \sum_{w\in W_3^g(s,t)}
    {1\over \vert{\rm Aut}(w)\vert}\quad ,\hfill \cr
    m_6^g(3r;2s,2t)\hfill&=&\displaystyle\sum_{w\in W_6^g(3r;2s,2t)}
    {1\over \vert{\rm Aut}(w)\vert}\quad .\hfill }$$

\par {\bf Theorem 1.1.\cite{[BV]}}
{\sl
\ \ 
(i) The group ${\rm Aut}(w)$ of automorphisms of an oriented
 maximal Wicks form $w$ is cyclic of order $1,\ 2,\ 3$ or $6$.
\par
 \ \ (ii) $\displaystyle
 m_1^g={2\over 1}\Big({1^2\over 12}\Big)^g{(6g-5)!\over g!(3g-3)!}
 \quad .$
\par \ \ (iii) $m_2^g(r)>0$ (with $r\in {\bf N}$) if and only if
 $f={2g+1-r\over 4}\in \{0,1,2,\dots\}$ and we have then
\par $\displaystyle m^g_2(r)={2\over 2}\Big({2^2\over 12}\Big)^f
 {1\over r!}{(6f+2r-5)!\over f!(3f+r-3)!}
 \quad .$
\par \ \ (iv) $m_3^g(s,t)>0$ if and only if
 $f={g+1-s-t\over 3}\in \{0,1,2,\dots\}$, $s\equiv 2g+1\pmod 3$ and
 $t\equiv 2g\pmod 3$ (which follows from the two previous conditions).
 We have then
\par $\displaystyle m^g_3(s,t)={2\over 3}\Big({3^2\over 12}\Big)^f
 {1\over s!t!}{(6f+2s+2t-5)!\over f!(3f+s+t-3)!}$ if $g>1$ and 
$\displaystyle m^1_3(0,2)={1\over 6}$.
\par \ \ (v) $m_6^g(3r;2s,2t)>0$ if and only if
 $f={2g+5-3r-4s-4t \over 12}\in \{0,1,2,\dots\}$,
 $2s\equiv 2g+1\pmod 3$ and $2t\equiv 2g\pmod 3$ (follows in fact
 from the previous conditions). We have then
\par $\displaystyle
 m_6^g(3r;2s,2t)={2\over 6}\Big({6^2\over 12}\Big)^f
 {1\over r!s!t!}{(6f+2r+2s+2t-5)!\over f!(3f+r+s+t-3)!}$\\
if $g>1$ and $\displaystyle m^1_6(3;0,2)={1\over 6}$.}


\bigbreak

\par {\bf Theorem 1.2.}
{\sl\ \ (i) The group ${\rm Aut}(S_g)$ of automorphisms of an extremal
 surface is cyclic of order $1,\ 2,\ 3$ or $6$ in every genus.
\par \ \ (ii) There is a bijection between the number of equivalence
 classes of oriented maximal genus $g$ Wicks forms
and extremal genus $g$ surfaces for $g \geq 4$.}

\medbreak

\par Set
 $$\matrix{m_2^g&=&\sum_{r\in {\bf N},\ (2g+1-r)/4\in {\bf N}\cup \{0\}}
 m_2^g(r)\quad ,\hfill \cr
 m_3^g&=&\sum_{s,t\in {\bf N},\ (g+1-s-t)/3\in {\bf N}\cup \{0\},\ 
 s\equiv 2g+1 \pmod 3} m_3^g(s,t)\quad ,\hfill \cr
 m_6^g&=&\sum_{r,s,t\in {\bf N},\ (2g+5-3r-4s-4t)/12\in
 {\bf N}\cup \{0\},\ 2s \equiv 2g+1\pmod 3} m_6^g(3r;2s,2t)\hfill}$$
 (all sums are finite) and denote by $M_d^g$ the number of isometry
 classes of extremal genus $g$ surfaces having an
  automorphism of order $d$ (i.e. an automorphism group with order
 divisible by $d$).

\medbreak

\par {\bf Theorem 1.3.}
 {\sl We have
 $$\matrix{ M_1^g=m_1^g+m_2^g+2m_3^g+2m_6^g\quad ,\hfill \cr
 M_2^g=2m_2^g+4m_6^g\quad ,\hfill \cr
 M_3^g=3m_3^g+3m_6^g\quad ,\hfill \cr
 M_6^g=6m_6^g \hfill} $$
 and $M_d^g=0$ if $d$ is not a divisor of $6$.}

\medbreak

\par The number $M_1^g$ of this Theorem is  the number of
 extremal surfaces
 of genus $g$ for $g \geq 4$. The first 15 values $M_1^1,\dots,M_1^{15}$ are
displayed in the Table at the end of this paper.

\par The following result is an immediate consequence of Theorem 1.3.

\medbreak

\par {\bf Corollary 1.1.}
 {\sl For $g \geq 4$ there are exactly
\par $M_6^g$ nonisometric extremal surfaces with $6$
 automorphisms,\par $M_3^g-M_6^g$ nonisometric extremal surfaces with $3$
 automorphisms,\par $M_2^g-M_6^g$ nonisometric extremal surfaces with $2$
 automorphisms and\par $M_1^g-M_2^g-M_3^g+M_6^g$ nonisometric extremal
surfaces without non-trivial automorphisms.}

\medbreak

\par {\bf Remark.}  Computing masses amounts to enumerating
 pointed objects, i.e. linear words instead of cyclic words
in Definition 1.1. Their number is $(12g-6)m_d^g$, where
$d$ is $ 1,2,3$ or $6$.

\par Let us remark that formula (ii) can be obtained from
 \cite{[WL]} (formula (9) on page 207 and the formula on the top of page
 211) or from \cite{[GS]} (Theorem 2.1 with $\lambda=2^{6g-3}$
 and $\mu=3^{4g-2}$).  
 Related objects have also been considered in \cite{[HZ]}.


\section{Oriented Wicks forms}
\label{sec:2}


To understand better the structure of extremal surfaces,
we need to describe now some properties of Wicks forms
(see \cite{[V]}, \cite{[BV]}
for all the details).

\par Let $V$ be a negative vertex of an oriented maximal Wicks
form of genus $g>1$. There are three possibilities, denoted
configurations of type $\alpha,\ \beta$ and $\gamma$
(see Figure 1) for the local configuration at $V$.

\bigbreak

\centerline{\epsfysize2.5cm\epsfbox{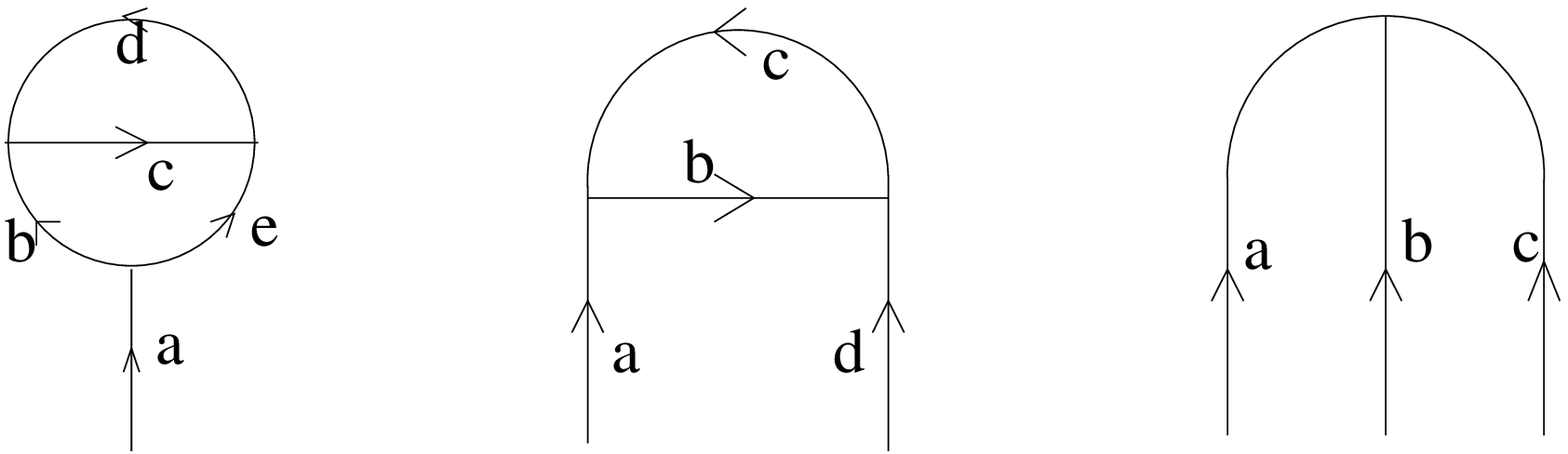}}
\centerline{\it Figure 1.}

\bigbreak

\par Type $\alpha$. The vertex $V$ has only two neighbours 
which are adjacent to each other. This implies that $w$ is of the
form
$$w=x_1abcdb^{-1}ec^{-1}d^{-1}e^{-1}a^{-1}x_2u_1x_2^{-1}x_1^{-1}u_2$$
(where $u_1,u_2$ are subfactors of $w$)
and $w$ is obtained from the maximal oriented Wicks form 
$$w'=xu_1x^{-1}u_2$$ 
of genus $g-1$ by the substitution  $x\longmapsto
x_1abcdb^{-1}ec^{-1}d^{-1}e^{-1}a^{-1}x_2$ and $x^{-1}\longmapsto x_2^{-1}
x_1^{-1}$ (this construction is called the $\alpha -$construction in
\cite{[V]}).

\par Type $\beta$. The vertex $V$ has two non-adjacent neighbours. 
The word $w$ is then of the form
$$w=x_1abca^{-1}x_2u_1y_1db^{-1}c^{-1}d^{-1}y_2u_2$$
(where perhaps $x_2=y_1$ or $x_1=y_2$, see \cite{[V]} for all the details).
The word $w$ is then obtained by a $\beta -$construction from the 
word $w'=xu_1yu_2$ which is an oriented maximal Wicks form 
of genus $g-1$.

\par Type $\gamma$. The vertex $V$ has three distinct neighbours. We have
then
$$w=x_1ab^{-1}y_2u_1z_1ca^{-1}x_2u_2y_1bc^{-1}z_2u_3$$
(some identifications among $x_i,\ y_j$ and $z_k$ may occur,
see \cite{[V]}
for all the details) and the word $w$ is obtained by a 
so-called $\gamma-$construction from the word $w'=x\tilde u_2y\tilde
u_1z\tilde u_3$.

\medbreak

\par {\bf Definition 2.1.}
We call the application which associates to an oriented maximal Wicks
form $w$ of genus $g$ with a chosen negative vertex $V$ the oriented
maximal Wicks form $w'$ of genus $g-1$ defined as above the {\it reduction}
of $w$ with respect to the negative vertex $V$. 

\medbreak

\par An inspection of figure $1$ shows that reductions with respect to
vertices of type $\alpha$ or $\beta$ are always paired since two doubly
adjacent vertices are negative, of the same type 
($\alpha$ or $\beta$) and yield the same reductions.

\par The above constructions of type $\alpha,\ \beta$ and $\gamma$
can be used for a recursive construction of all 
oriented maximal Wicks forms of genus $g>1$. 

\medbreak

\par {\bf Definition 2.2.} Consider an oriented maximal Wicks form
$$w=w_1\dots w_{12g-6}.$$ To any
edge $x$ of $w$ we associate a transformation of $w$ called the
{\it IH-transformation on the edge $x$}.

\medbreak

Geometrically, an IH-transformation amounts to contracting the
edge $x$ of the graph $\Gamma\subset S$ representing the oriented
maximal Wicks form $w$. 
This creates a vertex of degree $4$ which can be split
in two different ways (preserving planarity of the graph on $S$)
into two adjacent vertices of degree $3$: 
The first way gives back the original Wicks form and the 
second way results in the IH-transformation. Graphically, an
IH-transformation amounts hence to replace a (deformed) letter I 
(a topological neighbourhood of the edge $x\in \Gamma\subset S$)
by a (deformed) letter H.

\par More formally, one considers the two subfactors $axb$ and $cx^{-1}d$ 
of the (cyclic) word $w$. Geometric considerations and Definition 1.1
show that $b\not=a^{-1},c\not=b^{-1},d\not=a^{-1},d\not=c^{-1}$ and
$(c,d)\not=(a^{-1},b^{-1})$. 

\par According to the remaining possibilities we consider now the following
transformation:

\par \ \ Type 1. $c\not=a^{-1}$ and $d\not=b^{-1}$. This implies that
$d^{-1}a^{-1}$ and $b^{-1}c^{-1}$ appear as subfactors in the cyclic word 
$w$. The IH-transformation on the edge $x$ is then defined by the
substitutions
$$\matrix{axb&\longmapsto &ab\cr
cx^{-1}d&\longmapsto &cd\cr
d^{-1}a^{-1}&\longmapsto&d^{-1}ya^{-1}\cr
b^{-1}c^{-1}&\longmapsto&b^{-1}y^{-1}c^{-1}}$$
in the cyclic word $w$.

\par \ \ Type 2a. Suppose $c^{-1}=a$. This implies that $b^{-1}axb$ and
$d^{-1}a^{-1}x^{-1}d$ are subfactors of the cyclic word $w$. Define 
the IH-transformation on the edge $x$ by
$$\matrix{b^{-1}axb&\longmapsto& b^{-1}yab\hfill\cr
d^{-1}a^{-1}x^{-1}d&\longmapsto&d^{-1}y^{-1}a^{-1}d\quad .}$$

\par \ \ Type 2b. Suppose $d^{-1}=b$. Then $axba^{-1}$ and $cx^{-1}b^{-1}
c^{-1}$ are subfactors of the cyclic word $w$ and we define the IH-
transformation on the edge $x$ by
$$\matrix{axba^{-1}&\longmapsto&abya^{-1}\hfill\cr
cx^{-1}b^{-1}c^{-1}&\longmapsto&cb^{-1}y^{-1}c^{-1}\quad .}$$

\medbreak

\par {\bf Lemma 2.1.} {\sl (i) 
IH-transformations preserve oriented maximal 
Wicks forms of genus $g$. 

\par (ii) Two oriented maximal Wicks forms related by a IH-transformation
of type 2 are equivalent.}

\par Proof. This results easily by considering the effect of
an IH-transformation on the graph $\Gamma\subset S$.\hfill QED

\medbreak

\par {\bf Proposition 2.1.} {\sl An oriented maximal Wicks 
form of genus $g$ has exactly $2(g-1)$ positive and $2g$ negative 
vertices.}

\medbreak

\par {\bf Lemma 2.2.} {\sl An $\alpha$ or a $\beta$ construction increases
the number of positive and negative vertices by 2.}

\par The proof is easy.

\medbreak

\par {\bf Lemma 2.3.} {\sl The number of positive or negative vertices is
constant under IH-transformations.}

\par Proof of Lemma 2.3. The Lemma holds for IH-transformations of type 2
by Lemma 2.1 (ii). Let hence $w,w'$ be two oriented maximal Wicks forms
related by an IH-transformation of type 1 with respect to the edge
$x$ of $w$ respectively $y$ of $w'$. This implies that $w$ contains the
four subfactors
$$axb\quad ,\quad cx^{-1}d\quad ,\quad d^{-1}a^{-1}\quad ,\quad b^{-1}
c^{-1}$$
and $w'$ contains the subfactors
$$ab\quad ,\quad cd\quad ,\quad d^{-1}ya^{-1}\quad ,\quad b^{-1}y^{-1}
c^{-1}$$
in the same cyclic order and they agree everywhere else. It is hence
enough to check the lemma for the six possible cyclic orders of the
above subfactors.

\par One case is
$$\matrix {
w\hfill
&=&axbu&\dots& cx^{-1}d&\dots& d^{-1}a^{-1}&\dots& b^{-1}c^{-1}&\dots\quad ,
\hfill \cr
w'\hfill
&=&abu&\dots& cd&\dots& d^{-1}ya^{-1}&\dots& b^{-1}y^{-1}c^{-1}&\dots\quad ,
\hfill\quad . }$$
In this case the two vertices of $w$ incident in $x$ and the two vertices
of $w'$ incident in $y$ have opposite signs. All other vertices are
not involved in the IH-transformation and keep their sign and the Lemma
holds hence in this case.

\par The five remaining cases are similar and left to the reader.
\hfill QED

\medskip
\par Proof of Proposition 2.1. The result is true in genus 1 by inspection
(the cyclic word $a_1a_2a_3a_1^{-1}a_2^{-1}a_3^{-1}$ is the unique 
oriented maximal Wicks form of genus 1 and has two negative vertices.) 

\par Consider now an oriented maximal Wicks form $w$ of genus $g+1$. 
Choose an oriented embedded loop $\lambda$
of minimal (combinatorial) length in $\Gamma$. 

\par First case.
If $\lambda$ is of length 2 there are two vertices related by a 
double edge in $\Gamma$. This implies that they are negative and of type
$\alpha$ or $\beta$. The assertion of Proposition 2.1 holds hence for 
$w$ by Lemma 2.2 and by induction on $g$.

\par Second case.
We suppose now that $\lambda$ is of length $\geq 3$. The oriented
loop $\lambda$ turns either left or right at each vertex.
If it turns on the same side at two consecutive vertices $V_i$ and $V_{i+1}$
the IH-transformation with respect to the edge joining $V_i$ and $V_{i+1}$
transforms $w$ into a form $w'$ containing a shorter loop. By Lemma
2.2, the oriented maximal Wicks forms $w$ and $w'$ have the same number 
of positive (respectively negative)  vertices.
\par 
If $\lambda$ does not contain two consecutive vertices $V_i$ and $V_{i+1}$ 
with the above 
property (ie. if $\lambda$ turns first left, then right, then left etc.)
choose any edge $\{V_i,V_{i+1}\}$ in $\lambda$ and make an 
IH-transformation with respect to
this edge. This produces a form $w'$ which contains a loop $\lambda'$ of
the same length as $\lambda$ but turning on the same side at the
two consecutive vertices $V_{i-1},V_i$ or $V_{i+1},V_{i+2}$. 
By induction on the length of $\lambda$ we can
hence relate $w$ by a sequence of IH-transformation to an oriented
maximal Wicks form $\tilde w$ of genus $g+1$ containing a loop of length 2 
for which the result holds by the first case. The Wicks forms $w$ and
$\tilde w$ have of course the same number of positive (respectively
negative) vertices by Lemma 2.2. \hfill QED

\section{Proof of Theorem 1.2}
\label{sec:3}

\par Proof of Theorem 1.2. 

\par Proof of (i).

\par Let $w$ be an oriented maximal Wicks form with an automorphism
$\mu$ of order $d$. Let $p$ be a prime dividing $d$. The automorphism
$\mu'=\mu^{d/p}$ is hence of order $p$. If $p\not= 3$ then $\mu'$
acts without fixed vertices on $w$ and Proposition 2.1 shows that $p$
divides the integers $2(g-1)$ and $2g$ which implies $p=2$. 
The order $d$ of $\mu$ is hence 
of the form $d=2^a3^b$. Repeating the above argument with the prime 
power $p=4$ shows
that $a\leq 1$. 

\par All orbits of $\mu^{2^a}$ on the set
of positive (respectively negative)
vertices have either $3^b$ or $3^{b-1}$ elements and this leads 
to a contradiction if $b\geq 2$. 
This shows that $d$ divides $6$ and proves that 
the automorphism groups of oriented maximal
 Wicks forms
are always cyclic of order 1,2,3 or 6.

Let's consider an extremal genus $g$ surface $S_g$.
It was proved in \cite{[Ba]}, that surface is extremal
if and only if it can be obtained from a regular hyperbolic
$12g-6$-gon  with angles $2\pi/3$ such that the image
of the boundary of the polygon after identification
of corresponded sides is a geodesic graph with $4g-2$ vertices
of valence $3$ and $6g-3$ edges of equal length.
It was shown in $\cite{[GG3]}$, that any isometry of
an extremal surface of genus $g>3$ is realized by
a rotation of the $12g-6$-gon.  

Let $P$ be a regular geodesic hyperbolic polygon with $12g-6$ equal sides
and all angles equal to $2\pi/3$, equipped with a oriented maximal
genus $g$ Wicks form $W$ on its boundary. Consider the surface $S_g$
obtained from $P$ by identification of sides with the same labels. 
Since we made the identification using an oriented maximal
 Wicks form of length
$12g-6$, then the boundary of $D$ becomes a graph $G$
 with $4g-2$ vertices
of valence $3$ and $6g-3$ edges(see section 2).  
We started from a regular geodesic hyperbolic polygon with angles
$2\pi/3$, so $G$ is a geodesic graph with edges of equal length.
By a result of C.Bavard, mentioned above, the surface
$S_g$ is extremal.

So, the surface is extremal if and only if it
can be obtained from a regular hyperbolic polygon with $12g-6$ equal sides
and all angles equal to $2\pi/3$, equipped with an oriented maximal
genus $g$ oriented maximal Wicks form $W$ on its boundary.
The isometry of $S_g$ must be realized by a rotation
of the $12g-6$-gon \cite{[GG3]}, so the isometry must be an automorphism
of the Wicks form. Since the automorphism groups of oriented maximal
 Wicks forms
are always cyclic of order 1,2,3 or 6 then the automorphism
groups of genus $g \geq 4$ extremal surfaces are also
 cyclic of order 1,2,3 or 6.

\par Proof of (ii).

Every oriented maximal Wicks form defines exactly one extremal
surface, namely the surface obtained from 
 a regular hyperbolic polygon with $12g-6$ equal sides
and all angles equal to $2\pi/3$ with an oriented maximal
genus $g$ oriented maximal Wicks form $W$ on its boundary.

So, to prove the bijection between the number of equivalence
 classes of oriented maximal genus $g$ Wicks forms
and extremal genus $g$ surfaces for $g \geq 4$ we need
to show, that for every extremal surface $S_g$ there is only
one oriented maximal Wicks $W$ form such that $S_g$
can be obtained from a regular hyperbolic $12g-6$-gon
with $W$ on its boundary. It was proved in \cite{[GG3]}, that
for $g \geq 4$ the extremal disk $D$ of radius 
 $R_g=cosh^{-1}(1/sin \beta_g)$, $\beta_g=\pi/{(12g-6)}$,
 embedded in $S_g$ is unique. 
Consider the center $c$ of the disk $D$. The discs of radius
$R_g$ with the centers in the images of $c$
in the universal covering of $S_g$ form a packing of
the hyperbolic plane by discs.
To this packing one can classicly associate a tesselation $T$
of the hyperbolic plane by regular $12-6$-gons, which
are Dirichlet domains for $S_g$. And such a tesselation
is unique because of negative curvature (see \cite{[Ba]}, for example).
But each $T$ defines a Wicks form of length $12g-6$
in a unique way.

\par Theorem 1.2 is proved. \hfill QED

\par {\bf Table.}
The number of extremal surfaces in genus $1,2,4 \dots,15$:
$$
\matrix{
1\hfill &  1\hfill \cr 
2\hfill &  9\hfill \cr 
4\hfill &  1349005\hfill \cr 
5\hfill &  2169056374\hfill \cr 
6\hfill &  5849686966988\hfill \cr 
7\hfill &  23808202021448662\hfill \cr 
8\hfill &  136415042681045401661\hfill \cr 
9\hfill &  1047212810636411989605202\hfill \cr 
10\hfill &  10378926166167927379808819918\hfill \cr 
11\hfill &  129040245485216017874985276329588\hfill \cr 
12\hfill &  1966895941808403901421322270340417352\hfill \cr 
13\hfill &  36072568973390464496963227953956789552404\hfill \cr 
14\hfill &  783676560946907841153290887110277871996495020\hfill \cr 
15\hfill &  19903817294929565349602352185144632327980494486370\hfill \cr 
}$$

\end{document}